\documentclass[12pt]{article}

\usepackage{amsmath, epsfig, cite}
\usepackage{amssymb}
\usepackage{amsfonts}
\usepackage{latexsym}
\usepackage{float}
\usepackage{color}

\newtheorem{thm}{Theorem}[section]

\newtheorem{lem}[thm]{Lemma}

\numberwithin{equation}{section}

\newcommand{\qed}{{\hfill$\square$}\medskip}

\begin{document}

\begin{center}
{\Large\bf A variation of $q$-Wolstenholme's theorem}
\end{center}

\vskip 2mm \centerline{Ji-Cai Liu}
\begin{center}
{\footnotesize Department of Mathematics, Wenzhou University, Wenzhou 325035, PR China\\
{\tt jcliu2016@gmail.com} }
\end{center}


\vskip 0.7cm \noindent{\bf Abstract.}
We investigate a variation of $q$-Wolstenholme's theorem, which extends the $q$-analogue of Wolstenholme's theorem due to Shi and Pan [Amer. Math. Monthly 114 (2007), 529--531]. The proof makes use of the Ramanujan sum and higher order logarithmic derivatives of cyclotomic polynomials.


\vskip 2mm
\noindent{\it MR Subject Classifications}: 11A07, 11B65, 11T06

\section{Introduction}
In 1862, Wolstenholme asserted that if $p\ge 5$ is a prime, then the numerator of the fraction
\begin{align*}
1+\frac{1}{2}+\frac{1}{3}+\cdots+\frac{1}{p-1}
\end{align*}
written in the reduced form is divisible by $p^2$. Nowadays this result is known as the famous
Wolstenholme's theorem, which plays an important role in number theory. We refer to \cite{mestrovic-a-2011} for recent development and various extensions of Wolstenholme's theorem.

This paper focuses on $q$-analogues of the Wolstenholme's theorem. To continue the $q$-story of
Wolstenholme's theorem, we need some necessary notation.
For polynomials $A_1(q), A_2(q),P(q)\in \mathbb{Z}[q]$, the $q$-congruence $$A_1(q)/A_2(q)\equiv 0\pmod{P(q)}$$ is understood as $A_1(q)$ is divisible by $P(q)$ and $A_2(q)$ is coprime with $P(q)$. In general, for rational functions $A(q),B(q)\in \mathbb{Z}(q)$,
\begin{align*}
A(q)\equiv B(q)\pmod{P(q)}\Longleftrightarrow
A(q)-B(q)\equiv 0\pmod{P(q)}.
\end{align*}

Throughout the paper, let $\zeta$ denote a primitive $n$th root of unity.
The $q$-integers are defined as $[n]=(1-q^n)/(1-q)$ for $n\ge 1$, and the $n$th
cyclotomic polynomial is given by
\begin{align}
\Phi_n(q)=\prod_{\substack{1\le k \le n\\[3pt](n,k)=1}}
(x-\zeta^k).\label{new-1}
\end{align}
It is obvious that $\Phi_p(q)=[p]$ for any prime $p\ge 3$ and
\begin{align}
A(q)\equiv 0\pmod{\Phi_n(q)}
\Longleftrightarrow A(\zeta)=0,\label{new-2}
\end{align}
for rational function $A(q)$.

In the past few years, $q$-analogues of congruences ($q$-congruence) for indefinite sums of binomial coefficients as well as hypergeometric series attracted many experts' attention
(see, for example, \cite{guo-aam-2020,gs-rm-2019,gs-ca-2019,gz-am-2019,lp-a-2019}).

In 2007, Shi and Pan \cite{sp-amm-2007} established an interesting $q$-analogue of Wolstenholme's theorem as follows:
\begin{align}
\sum_{k=1}^{p-1}\frac{1}{[k]}\equiv \frac{(p-1)(q-1)}{2}+\frac{(p^2-1)(1-q)(1-q^p)}{24}
\pmod{[p]^2},\label{new-3}
\end{align}
for any prime $p\ge 5$. Letting $q\to 1$ on both sides of \eqref{new-3}, we reach the
Wolstenholme's theorem. The proof of \eqref{new-3} in \cite{sp-amm-2007} makes use of the following essential formula:
\begin{align}
\sum_{k=1}^{p-1}\frac{\zeta^k}{(1-\zeta^k)^2}=\frac{1-p^2}{12}.\label{new-4}
\end{align}

The motivation of the paper is to give an extension of \eqref{new-4}, and utilize this generalized result to obtain another $q$-analogue of Wolstenholme's theorem, which extends \eqref{new-3}.

Before stating the main results, we recall some necessary arithmetic functions.
The M\"obius function $\mu(n)$ is an important multiplicative function in number theory, which is
given by
\begin{align*}
\mu(n)
=\begin{cases}
(-1)^r\quad&\text{if $n$ is the product of $r$ different primes,}\\[5pt]
0\quad&\text{otherwise.}
\end{cases}
\end{align*}
The Euler totient function $\varphi(n)$ is the number of positive integers not greater than $n$ that are coprime to $n$, and the Jordan totient function $J_k(n)$ is the number of $k$-tuples of positive integers all less than or equal to $n$ that form a coprime $(k + 1)$-tuple together with $n$. Note that Jordan totient function is a generalization of Euler totient function:
\begin{align}
J_1(n)=\varphi(n).\label{a-1}
\end{align}
It is well-known that (see \cite[page 91]{siva-b-1989})
\begin{align*}
J_k(n)=\sum_{d|n}\mu(n/d)d^k.
\end{align*}

The first aim of the paper is to extend \eqref{new-4} as follows.
\begin{thm}\label{t-1}
Let $n\ge 2$ be an integer and $\zeta$ be a primitive $n$th root of unity. Then
\begin{align}
\sum_{\substack{1\le k \le n\\[3pt](n,k)=1}}\frac{\zeta^k}{(1-\zeta^k)^2}=-\frac{J_2(n)}{12}.\label{a-2}
\end{align}
\end{thm}

The second purpose of the paper is to utilize \eqref{a-2} to obtain the following $q$-congruence,
which generalizes \eqref{new-3}.
\begin{thm}\label{t-2}
For any positive integer $n\ge 2$, we have
\begin{align}
\sum_{\substack{1\le k \le n\\[3pt](n,k)=1}}\frac{1}{[k]}
\equiv\frac{(1-q)\varphi(n)}{2}+\frac{(1-q)(1-q^n)J_2(n)}{24}\pmod{[n]\Phi_n(q)}.\label{a-3}
\end{align}
\end{thm}

The important ingredients in the proof of \eqref{a-2} are higher order logarithmic derivatives
of cyclotomic polynomials and the following Ramanujan sum:
\begin{align*}
r_j(n)=\sum_{\substack{1\le k \le n\\[3pt](n,k)=1}}\zeta^{kj},
\end{align*}
which can also be expressed as
\begin{align}
r_j(n)=\frac{\varphi (n)\mu\left(n/(n,j)\right)}
{\varphi(n/(n,j))}.\label{c-2}
\end{align}
One refers to \cite[Section 16.6]{hw-b-1979} for the proof of \eqref{c-2}.

We recall an important theorem due to Herrera-Poyatos and Moree \cite{hm-a-2018} and establish a polynomial formula in the next section. Theorems \ref{t-1} and \ref{t-2} are proved in
Section 3.

\section{Preliminary results}
In order to prove Theorems \ref{t-1} and \ref{t-2}, we require the following two preliminary results.
\begin{lem} (See \cite[Theorem 3.3]{hm-a-2018}.)
For integers $n\ge 2$ and $k\ge 1$, we have
\begin{align}
\left.\frac{\mathrm{d}^k \ln\Phi_n(z)}{\mathrm{d} z^k}\right|_{z=1}=\sum_{j=1}^k
\frac{B_j(1)s(k,j)}{j}J_j(n),\label{b-1}
\end{align}
where the Bernoulli polynomials and Stirling numbers are given by
\begin{align*}
\frac{xe^{tx}}{e^x-1}=\sum_{k=0}^{\infty}B_k(t)\frac{x^k}{k!}
\end{align*}
and
\begin{align*}
x(x-1)(x-2)\cdots (x-k+1)=\sum_{j=0}^ks(k,j)x^j,
\end{align*}
respectively.
\end{lem}

\begin{lem}
For any positive integer $n\ge 2$, we have
\begin{align}
\sum_{k=1}^n\frac{\mu(n/(n,k))}{\varphi(n/(n,k))}z^{k-1}=\frac{z^n-1}{\varphi(n)}
\cdot\frac{\mathrm{d} \ln\Phi_n(z)}{\mathrm{d} z}.\label{b-2}
\end{align}
\end{lem}

{\noindent\it Proof.}
By \eqref{new-1} and the binomial theorem, we have
\begin{align}
\frac{\mathrm{d} \ln\Phi_n(z)}{\mathrm{d} z}&=
\sum_{\substack{1\le k \le n\\[3pt](n,k)=1}}\frac{1}{z-\zeta^k}\notag\\[5pt]
&=-\sum_{\substack{1\le k \le n\\[3pt](n,k)=1}}\sum_{j=0}^{\infty}\zeta^{-k(j+1)}z^j\notag\\[5pt]
&=-\sum_{j=0}^{\infty}z^j\sum_{\substack{1\le k \le n\\[3pt](n,k)=1}}\zeta^{-k(j+1)}\notag\\[5pt]
&=-\sum_{j=0}^{\infty}z^j\sum_{\substack{1\le k \le n\\[3pt](n,k)=1}}\zeta^{k(j+1)}.\label{b-3}
\end{align}
It follows from \eqref{b-3} and \eqref{c-2} that
\begin{align*}
\frac{\mathrm{d} \ln\Phi_n(z)}{\mathrm{d} z}
&=-\varphi (n)\sum_{j=0}^{\infty}\frac{\mu\left(n/(n,j+1)\right)}
{\varphi(n/(n,j+1))}z^j\\[5pt]
&=-\varphi (n)\sum_{j=1}^{\infty}\frac{\mu\left(n/(n,j)\right)}
{\varphi(n/(n,j))}z^{j-1}.
\end{align*}
Thus,
\begin{align*}
\frac{z^n-1}{\varphi(n)}
\cdot\frac{\mathrm{d} \ln\Phi_n(z)}{\mathrm{d} z}
&=\sum_{j=1}^{\infty}\frac{\mu\left(n/(n,j)\right)}
{\varphi(n/(n,j))}z^{j-1}-\sum_{j=1}^{\infty}\frac{\mu\left(n/(n,j)\right)}
{\varphi(n/(n,j))}z^{n+j-1}\\[5pt]
&=\sum_{j=1}^{\infty}\frac{\mu\left(n/(n,j)\right)}
{\varphi(n/(n,j))}z^{j-1}-\sum_{j=1}^{\infty}\frac{\mu\left(n/(n,n+j)\right)}
{\varphi(n/(n,n+j))}z^{n+j-1}\\[5pt]
&=\sum_{j=1}^{n}\frac{\mu\left(n/(n,j)\right)}
{\varphi(n/(n,j))}z^{j-1},
\end{align*}
where we have utilized the fact $(n,j)=(n,n+j)$ in the second step.
\qed

\section{Proof of the main results}
{\noindent\it Proof of Theorem \ref{t-1}.}
Let
\begin{align*}
f_n(z)=\sum_{\substack{1\le k \le n\\[3pt](n,k)=1}}\frac{\zeta^k}{(1-z\zeta^k)^2}.
\end{align*}
By the binomial theorem, we have
\begin{align}
f_n(z)&=\sum_{\substack{1\le k \le n\\[3pt](n,k)=1}}\zeta^k\sum_{j=0}^{\infty}{-2\choose j}(-z\zeta^k)^j\notag\\
&=\sum_{j=0}^{\infty}(j+1)z^j\sum_{\substack{1\le k \le n\\[3pt](n,k)=1}}\zeta^{k(j+1)}\notag\\
&=\sum_{j=1}^{\infty}jz^{j-1}\sum_{\substack{1\le k \le n\\[3pt](n,k)=1}}\zeta^{kj}.\label{c-1}
\end{align}
Substituting \eqref{c-2} into \eqref{c-1} gives
\begin{align*}
f_n(z)=\varphi (n)\sum_{j=1}^{\infty}\frac{\mu\left(n/(n,j)\right)}
{\varphi(n/(n,j))}jz^{j-1}.
\end{align*}
It follows that
\begin{align}
f_n(z)
&=\varphi (n)\sum_{j=0}^{\infty}\sum_{s=1}^{n}\frac{\mu\left(n/(n,nj+s)\right)}
{\varphi(n/(n,nj+s))}(nj+s)z^{nj+s-1}\notag\\
&=\varphi (n)\left(\sum_{j=0}^{\infty}njz^{nj-1}\sum_{s=1}^{n}\frac{\mu\left(n/(n,s)\right)}
{\varphi(n/(n,s))}z^{s}+\sum_{j=0}^{\infty}z^{nj-1}\sum_{s=1}^{n}\frac{\mu\left(n/(n,s)\right)}
{\varphi(n/(n,s))}sz^{s}\right)\notag\\
&=\varphi (n)\left(\frac{nz^{n-1}}{(1-z^n)^2}\sum_{s=1}^{n}\frac{\mu\left(n/(n,s)\right)}
{\varphi(n/(n,s))}z^{s}+\frac{1}{z(1-z^n)}\sum_{s=1}^{n}\frac{\mu\left(n/(n,s)\right)}
{\varphi(n/(n,s))}sz^{s}\right),\label{c-3}
\end{align}
where we have utilized the fact $(n,nj+s)=(n,s)$ in the second step.

Let $L_n(z)$ denote the left-hand side of \eqref{b-2}:
\begin{align*}
L_n(z)=\sum_{k=1}^n\frac{\mu(n/(n,k))}{\varphi(n/(n,k))}z^{k-1}.
\end{align*}
It follows from \eqref{a-1}, \eqref{b-1} and \eqref{b-2} that
\begin{align}
&\left.L_n(z)\right|_{z=1}=0,\label{c-4}\\[10pt]
&\left.L_n^{'}(z)\right|_{z=1}=\frac{n}{2},\label{c-5}\\[10pt]
&\left.L_n^{''}(z)\right|_{z=1}=\frac{n(n-3)}{2}+\frac{nJ_2(n)}{6\varphi(n)}.\label{c-6}
\end{align}

Furthermore, we rewrite \eqref{c-3} as
\begin{align}
f_n(z)&=\varphi(n)\left(\frac{nz^nL_n(z)}{(1-z^n)^2}+\frac{(zL_n(z))^{'}}{1-z^n}\right)\notag\\[10pt]
&=\varphi(n)\cdot \frac{((n-1)z^n+1)L_n(z)+z(1-z^n)L_n^{'}(z)}{(1-z^n)^2}.\label{c-7}
\end{align}
By using the L'Hospital's rule and \eqref{c-4}--\eqref{c-7}, we obtain
\begin{align*}
\lim_{z\to 1}f_n(z)&=\varphi(n)\lim_{z\to 1} \frac{((n-1)z^n+1)L_n(z)+z(1-z^n)L_n^{'}(z)}{(1-z^n)^2}\\[10pt]
&=-\frac{\varphi(n)}{2n}\left(n(n-1)\lim_{z\to 1} \frac{L_n(z)}{1-z^n}+
2\left.L_n^{'}(z)\right|_{z=1}+\left.L_n^{''}(z)\right|_{z=1}\right)\\[10pt]
&=-\frac{\varphi(n)}{2n}\left((1-n)\left.L_n^{'}(z)\right|_{z=1}+
2\left.L_n^{'}(z)\right|_{z=1}+\left.L_n^{''}(z)\right|_{z=1}\right)\\[10pt]
&=-\frac{J_2(n)}{12},
\end{align*}
which implies
\begin{align*}
\sum_{\substack{1\le k \le n\\[3pt](n,k)=1}}\frac{\zeta^k}{(1-\zeta^k)^2}=-\frac{J_2(n)}{12},
\end{align*}
as desired.
\qed

{\noindent\it Proof of Theorem \ref{t-2}.}
For $1\le k\le n$, we have
\begin{align*}
(n,k)=1\Longleftrightarrow (n,n-k)=1.
\end{align*}
Thus,
\begin{align*}
\sum_{\substack{1\le k \le n\\[3pt](n,k)=1}}\frac{1}{[k]}
&=(1-q)\sum_{\substack{1\le k \le n\\[3pt](n,k)=1}}\frac{1}{1-q^k}\\[10pt]
&=\frac{1-q}{2}\sum_{\substack{1\le k \le n\\[3pt](n,k)=1}}\left(\frac{1}{1-q^k}+\frac{1}{1-q^{n-k}}\right)\\[10pt]
&=\frac{1-q}{2}\sum_{\substack{1\le k \le n\\[3pt](n,k)=1}}\left(\frac{1}{1-q^k}+\frac{1}{1-q^{n-k}}-1\right)+\frac{(1-q)\varphi(n)}{2}\\[10pt]
&=\frac{(1-q)(1-q^n)}{2}\sum_{\substack{1\le k \le n\\[3pt](n,k)=1}}\frac{1}{(1-q^k)(1-q^{n-k})}
+\frac{(1-q)\varphi(n)}{2}.
\end{align*}
Since $1-q^n\equiv 0\pmod{[n]}$, we have
\begin{align}
\sum_{\substack{1\le k \le n\\[3pt](n,k)=1}}\frac{1}{[k]}\equiv -\frac{(1-q)(1-q^n)}{2}\sum_{\substack{1\le k \le n\\[3pt](n,k)=1}}\frac{q^k}{(1-q^k)^2}
+\frac{(1-q)\varphi(n)}{2}\pmod{[n]^2}.\label{c-8}
\end{align}

From \eqref{new-2} and \eqref{a-2}, we deduce that
\begin{align}
\sum_{\substack{1\le k \le n\\[3pt](n,k)=1}}\frac{q^k}{(1-q^k)^2}\equiv -\frac{J_2(n)}{12}
\pmod{\Phi_n(q)}.\label{c-9}
\end{align}
It follows from \eqref{c-8} and \eqref{c-9} that
\begin{align*}
\sum_{\substack{1\le k \le n\\[3pt](n,k)=1}}\frac{1}{[k]}
\equiv\frac{(1-q)\varphi(n)}{2}+\frac{(1-q)(1-q^n)J_2(n)}{24}\pmod{[n]\Phi_n(q)},
\end{align*}
as desired.
\qed

\vskip 5mm \noindent{\bf Acknowledgments.}
This work was supported by the National Natural Science Foundation of China (grant 11801417).

\end{document}